# ON ITERATING OPERATORS AND ON GENERALIZED PERIODIC ORBITS

Andrei Vieru

**Abstract**

We try to define the more general form of iterative processes in which the Pomeau-Manneville and the Feigenbaum scenario may occur along with their specific scaling properties. Doing this we need to generalize other basic concepts. Thus, what we call a periodic carousel is a generalization of what is usually called a periodic orbit.

## 1. ITERATING OPERATORS

### 1.1. GENERALIZED ITERATION

In our paper 'Generalized iteration, catastrophes and generalized Sharkovsky's ordering' we introduced the concepts of generalized iteration of the first and of the second kind of a function that maps a topological space $T^n$ onto a topological space $T^m$ ($n \geq m$). Generalized iteration of the second kind may be viewed as a generalization of usual iteration because it coincides with usual iteration when $m = n$. Although this attempt to generalize iteration leads to interesting results – some of whom are analogous to classical results in the field of Dynamical Systems – it seems not fully satisfactory inasmuch, on one hand, we did not yet propose a generalization which should include the already proposed ones as special cases, and inasmuch, on the other hand, it may conduct to the wrong idea that iteration of a map $f : T \rightarrow T$ (where T is a one-dimensional topological space) cannot be generalized. Besides, generalized iteration of the first kind can be further generalized without losing scaling properties we shall talk about in this paper.

Therefore we shall propose another generalization of the operator called *iteration*. We'll define the concept of *iterating operator*, which embraces both iteration of the first and second kind (see the aforementioned article) and, as well, recursive processes analogous to those studied in our article 'Short note on additive sequences and recursive processes'. In particular we'll see that iterating operators may work with functions of one single real argument.

## 1.2 ITERATING OPERATORS. A GENERAL DEFINITION

### 1.2.1. PRELIMINARY CONVENTIONS

We shall start this section with the simplest example: the example of mere iteration itself. Let $y = f(x)$ be some function that maps a set A onto itself. If we define the operator Ω (acting on the set of functions that map A onto itself) as $\Omega(f) = f \circ f$ and the iterates of Ω as $(\Omega^n(f))(x) = f^{n+1}(x)$, then the operator Ω is, in fact, iteration itself. Considering one of the usual ways of writing the iterates of $f(x)$ – i.e. $f(f(x))$, $f(f(f(x)))$ – we see that from a linguistic point of view at each step $x$ is replaced by $f(x)$.

In the general definition of an *iterating operator* we shall suppose, as a convention, that functions are written correctly without blanks, in particular without blanks after comas: $f(x,y)$ is correct, while $f(x,\quad y)$ isn't. We'll consider that one has already rigorously established what exactly is a correct written function formula (in the sense of formal grammars or formal systems).

### 1.2.2. A GENERAL DEFINITION

Let **B** be a set of symbols needed to correctly write down a function that maps the topological space $T^n$ onto a topological space $T^m$ $(m \leq n)$.

An iterating operator Ω is a set of formal grammar rules that transform correctly written formulae into *longer* correctly written formulae. When we iterate the action of these rules they are supposed to act simultaneously wherever they can in the formula they act on. **In other words an iterating operator is a Lindenmayer System whose strings are interpreted as mappings of a metric space onto itself or, into another metric space of lower or at most equal dimension.** We shall always consider that $\Omega^0(f) = f$

### 1.2.3. EXAMPLES

**Example 1**

Let $f$ map a topological space $T^2$ into T and let $\Omega$ transform $f(x,y)$ in $f(y,f(x,y))$. via the two grammar rules: $x \to y$ and $y \to f(x,y)$

So, the next iterate will be $(\Omega^2(f))(x,y) = f(f(x,y), f(y, f(x, y)))$

Fibonacci numbers are constructed using exactly these two grammar rules with the convention $f(x,y) = x + y$.

Generalized iteration might be recalled of as a special case of an iterating operator.

**Example 2**

Let $G(x,y,z)$ and $H(x,y,z)$ be maps from $T^3$ to T.

Let $F(x,y,z) = (G(x,y,z),H(x,y,z))$ be a map from $T^3$ to $T^2$.

Let the *iterating operator* $\Omega$ be the set of the five following context-depending[1] formal grammar rules:

$(y, \to (G(x,y,z),$

$,z, \to ,H(x,y,z),$

$(z, \to (H(x,y,z),$

$G(x,y,z) \to G(y,z,G(x, y, z))$

$H(x,y,z) \to H(z,G(x,y,z),H(x,y,z)$

---

[1] for example, following the first rule, $y$ is replaced by $G(x,y,z)$ whenever it is placed between an opening parenthesis and a coma

Then one can write, for example, the first three iterates of $\Omega$ as:

$$\Omega(F(x,y,z)) \equiv \Omega(G(x,y,z),H(x,y,z)) = (G(y,z,G(x,y,z)),H(z,G(x,y,z),H(x,y,z)))$$

$$\Omega^2(F(x,y,z))=$$
$$(G(G(x,y,z),H(x,y,z),G(y,z,G(x,y,z))),H(H(x,y,z),G(y,z,G(x,y,z)),H(z,G(x,y,z),H(x,y,z))))$$

$$\Omega^3(F(x,y,z)) =$$
$$(\mathrm{G}(G(y,z,G(x,y,z)),H(z,G(x,y,z),H(x,y,z)),G(G(x,y,z),H(x,y,z),G(y,z,G(x,y,z)))),H(H(z,G(x,y,z),H(x,y,z)),G(G(x,y,z),H(x,y,z),G(y,z,G(x,y,z))),H(H(x,y,z),G(y,z,G(x,y,z)),H(z,G(x,y,z),H(x,y,z)))))$$

(And so on… Note that this is an example of what we have formerly called Generalized iteration of the first kind.)

**Example 3**

Let $F(x, z)$ be a map from $\mathrm{T}^2$ to T and let $(\Omega(F))(x,z) = F(y,F(x,z))$

Let the iterating operator $\Omega$ be defined by the transformational grammar rules: $x \rightarrow y$, $y \rightarrow z$ and $z \rightarrow F(x,z)$. Adding the rule we get

$(\Omega^2(F))(x,z) = F(z,F(y,F(x,z)))$

$(\Omega^3(F))(x,z) = F(F(x,z),F(z,F(y,F(x,z))))$

etc.

**Example 4**

Let $F(x,y)$ be a map from $\mathrm{T}^2$ to T and let the iterating operator $\Omega$ be defined by the three context-free grammar rules:

$x \rightarrow y,\ \ y \rightarrow z,\ z \rightarrow F(x,y)$

Here the first iterates of the iterating operator $\Omega$ could be written as:

$(\Omega(F))(x,y) = F(y,z)$

$(\Omega^2(F))(x,y) = F(z,F(x,y))$

$(\Omega^3(F))(x,y) = F(F(x,y),F(y,z))$

$(\Omega^4(F))(x,y) = F(F(y\ z),F(z,F(x,y))$

**Example 5**

Let $G(x,y,z)$ and $H(x,y,z)$ be a maps from $T^2$ to T. Let $F(x,y,z) = (G(x,y,z),H(x,y,z))$ be a map from $T^3$ to $T^2$ and let the iterating operator Ω be defined by the three context-depending grammar rules:

$(z, \rightarrow (H(x,y,z),$

$G(x,y,z) \rightarrow G(\mathrm{z},G(x,y,z),H(x,y,z))$

$H(x,y,z) \rightarrow H(\mathrm{z},G(x,y,z),H(x,y,z))$

These three rules correspond (for the case $F\colon T^3 \rightarrow T^2$) to what we formerly called generalized iteration of the second kind

## 2. PERIODIC CAROUSELS.

### 2.1. MOTIVATION

The usual iteration operator usually gives birth to stable or repelling periodic orbits, whose elements are, naturally, supposed to be distinct.

The classical definition of a periodic orbit may be formulated as follows:

A *periodic orbit* with period $k$ for a map $g$

$x_{i+1} = g(x_i),\ x \in \boldsymbol{R}^n,\ n \geq 1$

is the set of $k$ distinct points $\{p_j = g^j(p_0) | j = 0,...,\ k-1\}$ with $\mathrm{g}^k(p_0) = p_0$ (Guckenheimer and Holmes, 1983)

Computations show that iterating operators also may engender cycles – stable or repelling. However, their elements ***are not necessarily distinct***. Therefore, a broader concept and a more general definition need to be introduced. We propose the following ones:

## 2.2. PERIODIC CAROUSELS. A DEFINITION

Let $\Omega$ be an iterating operator that transforms a map $g_i : \mathbf{T}^n \rightarrow \mathbf{T}^m$ $(m \leq n)$ into a map $\Omega(g_i)$: $\mathbf{T}^n \rightarrow \mathbf{T}^m$, $(m \geq n)$

For a map $g_0$ a **periodic carousel** with **least** period $k$ is the finite sequence of $k$ (not necessarily distinct) points in $\mathbf{T}^m$ $\{p_j\}_{0 \leq j \leq k-1}$ with $p_j = \left(\Omega^{j-1}(g_0)\right)\left(p_0\right)$ $(0 < j < k)$, with $\left(\Omega^{k-1}(g_0)\right)\left(p_0\right) = p_0$, with, for any non-negative integers $l$ and $p$,

$l \equiv p \pmod{k} \Rightarrow \left(\Omega^{l}(g_0)\right)\left(p_0\right) = \left(\Omega^{p}(g_0)\right)\left(p_0\right)$

and such that for any $k' < k$ there is at least a pair of distinct integers $l$ and $p$ with $l \equiv p \pmod{k'}$ and $\left(\Omega^{l}(g_0)\right)\left(p_0\right) \neq \left(\Omega^{p}(g_0)\right)\left(p_0\right)$

## 3. SCALING PROPERTIES IN A SLIGHTLY MODIFIED FEIGENBAUM SCENARIO

### 3.1 UNUSUAL BIFURCATIONS

As we showed in our papers Generalized Iteration and Short Note on Recursive Processes, dynamical systems based on iterating operators (with control parameters) display a doubling-period cascade before reaching Chaos. However, often – if not always – the first step in this process ***is not a usual bifurcation***.

To show only a few examples, with the increment of the control parameters (i.e. increasing at least one of them) the initial stable fixed point may split into a stable periodic carousel with five distinct points and with period five as well: take for instance the iterating operator $\Omega$ as it was defined in the first chapter of this paper, example 4, with the function $F_a(x,y) = ax(1-x)y(1-y)$ depending on the parameter $a$.

Other examples were provided[1] in which the initial fixed point splits into a n attracting periodic carousel with period 8 but with only 5 distinct elements: take for instance the iterating operator $\Omega$ as it was defined in the first chapter, example 3, with the function $F_a(x,y) = ax(1-x)y(1-y)$ depending on the parameter $a$.

---

[1] see our paper "Short note on additive sequences and on recursive processes" arxiv:0806.3972

In another example the initial fixed point splits into a periodic carousel with period 3 but with only 2 distinct elements: take for instance the iterating operator Ω as it was defined in the first chapter, example 1, with the function

$F_a(x,y) = ax(1–x)y(1–y)$ depending on the parameter $a$.

## 3.2. FEIGENBAUM'S SCALING PROPERTY SEEMS TO REMAIN UNCHANGED IN THE CONTEXT OF OPERATING OPERATORS

However, after the first somewhat unusual type of split, the increment of the control parameter yields a period-doubling cascade that leads to a chaos point, through a series of bifurcation points that converge to it at a rate that seems to coincide with the Feigenbaum δ constant for a large class of maps, regardless of the type of the chosen iterating operator (but under reasonable conditions, analogous — *mutatis mutandis* — to the condition of negative Schwarzian derivative in maps of the interval onto itself).

## 4. THE CASE OF ONE-DIMENSIONAL TOPOLOGICAL SPACES

One can imagine that the only one conceivable iterating operator that can act on maps with one single real variable is *iteration* itself. This idea is as wrong as the idea that a periodic helix is nothing more than a periodic orbit engendered by the mere iteration of a periodic function mod 1.

We shall reject these ideas showing simple examples of iterating operators other than mere iteration.

Let $f(x) = \text{frac}(0.4\sin(2\pi x)+x+a)$ (where $\text{frac}(x) = x - \lfloor x \rfloor$)

We have, for example, the following 'usual' iterates:

$f^2(x) = \text{frac}(0.4\sin(2\pi\ \text{frac}(0.4\sin(2\pi x)+x+a))+ \text{frac}(0.4\sin(2\pi x) + x + a) +a)$

$$f^3(x)=\mathrm{frac}(0.4\sin(2\pi(\mathrm{frac}(0.4\sin(2\pi(\mathrm{frac}(0.4\sin(2\pi x) + x + a))) + \mathrm{frac}(0.4\sin(2\pi x) + x + a) + a))) + \mathrm{frac}(0.4\sin(2\pi(\mathrm{frac}(0.4\sin(2\pi x) + x + a))) + \mathrm{frac}(0.4\sin(2\pi x) + x + a) + a) + a)$$

etc.

Here we deal with usual iteration of a circle map.

Let now $f(x) = \mathrm{frac}(0.4\sin(\pi x)+x+a)$

If we define $\Omega$ by the one single rule $x \rightarrow 0.4\sin(\pi x) + x + a$ then we'll have:

$(\Omega(f))(x) = \mathrm{frac}(0.4\sin(\pi(0.4\sin(\pi x) + x + a)) + x + a)$

$$(\Omega^2(f))(x) = \mathrm{frac}(0.4\sin(\pi(0.4\sin(\pi(0.4\sin(\pi x) + x + a)) + (0.4\sin(\pi x) + x + a) + a)) + (0.4\sin(\pi x) + x + a) + a)$$

etc.

Geometrical-topological interpretation: the operator $\Omega$ iterates a map from the circle onto a degenerate projective plan (a circle whose antipodal points are identified).

Another example liable to be interpreted in the same way:

Let $f_{A,B}(x) = A\sin^2(\sin 0.5\pi x)+x+B \pmod 2$ be a family of maps of the circle

(see picture below)

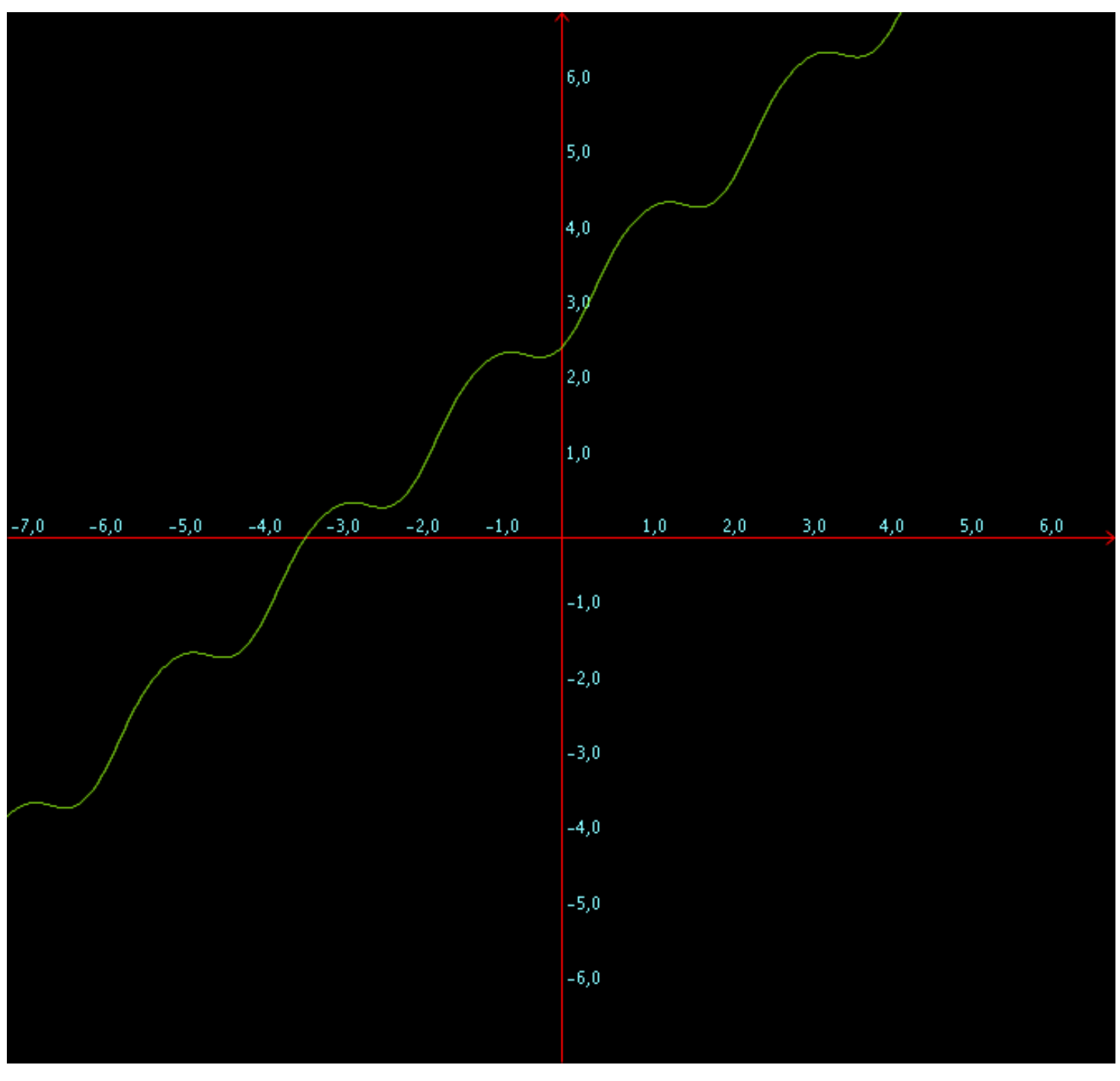

Let $(\Omega^n(f_{A,B}))(x) = \text{frac}(f^n{}_{A,B}(x))$

$(\Omega^n(f_{A,B}))(x)$ map the circle into a degenerate projective plane, just as in the previous example. Now let A=1.2 and B = 0.87907470162143…

Computations lead to a *stable periodic carousel* with period 3 and ***with only 2 distinct points*** (0.4595359333…, 0. 790594694…, 0.79059469…).

Remark: $f_{A,B}$ with fixed A and B values maps a one-dimensional topological space into another one-dimensional topological space. The iterating operator Ω acting on $f_{A,B}$ produces cycles: periodic orbits and periodic carousels. The question of their coexistence might deserve attention. ***Does a Sharkovsky-type ordering arise or not?***

## 5. SCALING PROPERTIES RELATED TO INTERMITTENCY AND TO THE POMEAU-MANNEVILLE SCENARIO

### 5.1. ITERATING OPERATORS, INTERMITTENCY AND POMEAU-MANNEVILLE SCENARIO

As one can naturally expect the behavior of most dynamical systems based on iterating operators applied to maps from the circle into a degenerate projective plan usually display intermittency before reaching full-scaled turbulent regime.

During computer simulations that finally led to a well-known result – only in a more general frame – we used the concept of *quasi-periodic pseudo-helix* with order *p*, length *m* and modulo 1 defined as follows. Basically, it corresponds to what is usually called 'laminar regime'.

Let *p* and *m* be integers. Let K=$\{\text{frac}(u_a(n_0)), \text{frac}(u_a(n_0+1)),\ldots, \text{frac}(u_a(n_0+m))\}$ be a finite sub-sequence of consecutive terms of some infinite sequence $\{\text{frac}(u_a(n))\}_{n\in\mathbb{N}}$.

We'll call K a *quasi-periodic pseudo-helix of order p with length m and with start s* if:

**A°)** the *p* finite sub-sub-sequences (functions whose variable *k* runs through a finite set of consecutive integer values)

$\Delta_{2,\,0,\,m,\,s}\,(k) = \{[\mathrm{frac}(u_a(s + p(k + 2))) - \mathrm{frac}(u_a(s + p(k + 1)))] - [\mathrm{frac}(u_a(s + p(k + 1))) -$ $\mathrm{frac}(u_a(s + pk))]\}_{0\,\le\,p\,(k+2)\,<\,m}$, **(1°)**

$\Delta_{2,\,1,\,m,\,s}\,(k)=\{[\mathrm{frac}(u_a(s+1+ p(k + 2))) - \mathrm{frac}(u_a(s + 1+ p(k + 1)))] - [\mathrm{frac}(u_a(s+1+ p(k +$ $1))) - \mathrm{frac}(u_a(s + 1 + pk))]\}_{0\,\le\,p\,(k+2)+\,1\,<\,m}$ **(2°)**

…………………

$\Delta_{2,\,p-1,\,m,\,s}\,(k) = \{[\mathrm{frac}(u_a(s + p - 1+ p(k + 2))) - \mathrm{frac}(u_a(s + p - 1+ p(k + 1)))] -$ $[\mathrm{frac}(u_a(s+ p - 1 + p(k + 1))) - \mathrm{frac}(u_a(s + p - 1+ pk))]\}_{0\,\le\,p\,(k+2)+\,p\,-\,1\,<\,m}$ **($p$°)**

are all strictly monotone.

**B°)** none of the infinite sub-sequences

$\Delta_{2,\,l,\,\infty,\,s}\,(\,k\,) =\{[\mathrm{frac}(u_a(s + l + p(k + 2))) - \mathrm{frac}(u_a(s+ l +p(k+1)))] - [\mathrm{frac}(u_a(s + l +p(k$ $+ 1))) - \mathrm{frac}(u_a(s + l + pk))]\}_{0\,\le\,p\,(k+2)+\,l\,<\,\infty,\;0\,\le\,l\,<\,p}$ converge to any limit.

**C°)** $\exists k_0 < m$ such that $\forall h \le k_0$ $\Delta_{2,\,0,\,m,\,s}(\,h\,)$, $\Delta_{2,\,1,\,m,\,s}\,(\,h\,)$,…, $\Delta_{2,\,p-1,\,m,\,s}\,(\,h\,)$ are of the same sign, while for any $i > k_0$ $\Delta_{2,\,0,\,m,\,s}(\,i\,)$, $\Delta_{2,\,1,\,m,\,s}\,(\,i\,)$,…, $\Delta_{2,\,p-1,\,m,\,s}\,(\,i\,)$ are also of the same (changed) sign.

**D°)** at least one of the finite sequences $\Delta_{2,\,0,\,m+1,\,s}\,(k)$, $\Delta_{2,\,0,\,m+1,\,s}\,(k\,)$,…, $\Delta_{2,\,0,\,m+1,\,s}\,(k)$ is not anymore monotone.

The computations are sometimes tricky so we shall not go through all their details here.

## 5.2. A SCALING PROPERTY RELATED TO INTERMITTENCY IN DYNAMICAL SYSTEMS GENERATED BY FUNCTIONS THAT MAP THE CIRCLE INTO A DEGENERATIVE PROJECTIVE PLAN

Computations show that

$$\lim_{\substack{P\to\infty \\ \lambda\to\infty}} \frac{\mu_F^{-1}(\lambda P)-\mu_F^{-1}(P)}{\mu_F^{-1}(\lambda^2 P)-\mu_F^{-1}(\lambda P)} = \lambda^2 \qquad \textbf{(1)}$$

where, for a given family of functions indexed on a control parameter, $\mu^{-1}{}_F$ is a supposedly one-to-one function that returns the value of the control parameter for which some average periodicity $\boldsymbol{P}$ of laminar phases is attained. Formula (1) is a possible expression of the frequency of laminar phases near the threshold were intermittency appears.

This numerically obtained result meets the classical result concerning quadratic scaling properties in usual dynamical systems based on iteration of circle maps (see, for example S.P. Kuznetsov, Dynamical Chaos, Chapter 17)

## 6. BASINS OF ATTRACTION

### 6.1.DYNAMICAL SYSTEMS BASED ON A MAPPING OF THE HYPERCUBE INTO AN INTERVAL

Let us consider a very simple example of a dynamical system based on the family of functions $f_a(x,y)=a\sin(\pi x)\sin(\pi x)$ and on the rules $x \to y$ and $y \to f(x,y)$

When the control parameter $a < 0.75$ we have a non zero stable fixed point whose basin of attraction looks like:

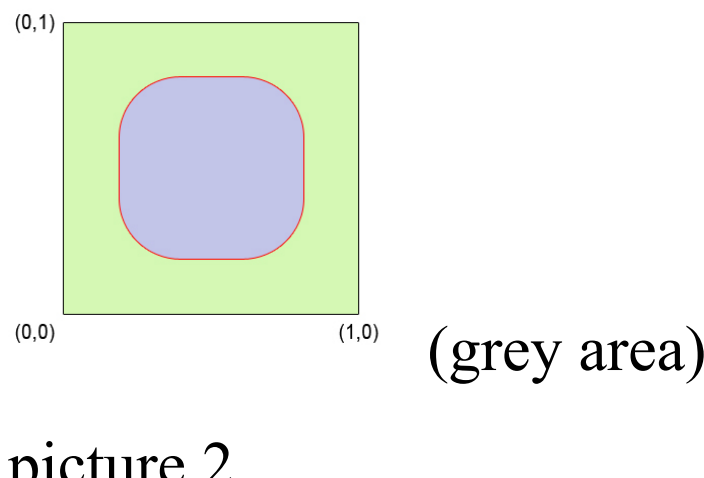

(grey area)

picture 2

When $a$ reaches (roughly) 0.754 we have a periodic carousel (0.7358383532881708…, 0.5496022243041646…, 0.5496022243041646…) whose basin of attraction is (probably) a fractal (see our paper 'Generalized iteration, catastrophes, generalized Sharkovsky ordering' where the considered family was $F_a(x,y)=ax(1-x)y(1-y)$.)

At $a\approx 0.765$ the basin of attraction is not anymore a fractal but also more or less of the shape of the blue area in the picture 2. At $a\approx 0.8102$ we have a bifurcation and a stable periodic carousel arises, with period 6 and with 6 distinct elements. Then chaos is quickly reached through a Feigenbaum-type process, where the Feigenbaum $\delta$ constant shows up. While the control parameter increases (after chaos point is already reached), we encounter an alternation of chaos regime with stable periodic carousels with period always divisible by 3.

When $0.90897 < a < 0.909057883$ we find a period 8 stable periodic carousel ($b_0$, $b_1$,$b_2$, $b_3$, $b_4$, $b_5$, $b_6$, $b_7$,) with tiny basins of attraction around points of the form ($b_i$, $b_{i+1}$)

(where the indexes are considered mod 8). For the rest of the initial values inside an area of the shape of the blue region in the picture 1, it seems that we have chaos with both a strange attractor and some kind of micro (or pseudo) intermittency.

### 6.2. BASINS OF ATTRACTION AND INTERMITTENCY WITH MAPS OF THE CIRCLE INTO A DEGENERATE PROJECTIVE PLAN

#### 6.2.1 FINITE TIME INTERMITTENCY (PSEUDO-INTERMITTENCY)

Usually intermittency occurs when in a circle map dissipative dynamical system the control parameter is set near some threshold between a laminar regime and a turbulent one. What happens on the threshold between two different laminar regimes? We shall provide an example of such an adjacent (on the scale of the control parameter) different laminar regimes.

Let the dynamical system be defined by $F_a(x)=$ frac($0.5\sin(\sin(\pi x))+x+a$) and by the rule $x \rightarrow 0.5\sin(\sin(\pi x))+x+a$

When $a = 0.8709$ we have a periodic helix with period 5

When $a = 0.870931$ we have a periodic helix with period 2

What happens between 0.8709 and 0.870931? Well, it depends on the initial condition. Let $a = 0.870931$ and the initial $x$ be 0.555. Then until about the 18175-th term, we'll have something that looks like intermittency with a quasi-periodic pseudo-helix of order 2, then a helix with period 5.

Let $a = 0.8709302$. For the randomly chosen values $x = 0.9$, $x = 0.611$, $x = 0.61$, $x = 0.3$, $x = 0.322$, $x = 0.3215$, $x = 0.3211$, $x = 0.32125$, $x = 0.321255$, $x = 0.329$, $x = 0.327$, $x = 0.326$, $x = 0.3265$, $x = 0.2$, $x = 0.22$, $x = 0.221$ we'll have helixes with period 2.

For the randomly chosen values $x = 0.601$, $x = 0.6005$, $x = 0.6$, $x = 0.5$, $x = 0.35$, $x = 0.34$, $x = 0.3333333333$, $x = 0.33$, $x = 0.328$, $x = 0.325$, $x = 0.32475$, $x = 0.32425$, $x = 0.3245$, $x = 0.324$, $x = 0.3231$, $x = 0.323$, $x = 0.321$, $x = 0.3212$, $x = 0.3213$, $x = 0.21$, $x = 0.205$, $x = 0.2075$, $x = 0.23$ we'll have helixes with period 5.

The structure of the set of initial values for which a helix with period 2 (respectively 5) arises remains unclear.

Another unanswered question is whether, playing with the value of the control parameter and/or with the initial condition the initial phase of (pseudo)-intermittency can be arbitrarily long.

## 7. A CONJECTURE ABOUT THE ROLE OF THE LINDENMAYER SYSTEMS IN THE PRESERVING OF THE CLASSICAL SCALING PROPERTIES

Further generalizations of the concept of iteration are possible (using formal grammar rules that are not Lindenmayer systems) but scaling properties related to the Feigenbaum or Pomeau-Manneville scenario will never arise.

andreivieru007@gmail.com

**REFERENCES**

(1) S.P. Kuznetsov, Dynamical Chaos, Chapter 17 (pages 257-259 in the Russian edition)

(2) Yves Pomeau and Paul Manneville, "Intermittent Transition to Turbulence in Dissipative Dynamical Systems", Commun. Math. Phys. vol. 74, pp. 189–197 1980

(3) M. J. Feigenbaum, "The Universal Metric Properties of Nonlinear Transformations." *J. Stat. Phys.* **21**, 669-706, 1979.

(4) M. J. Feigenbaum, "The Metric Universal Properties of Period Doubling Bifurcations and the Spectrum for a Route to Turbulence." *Ann. New York. Acad. Sci.* **357**, 330-336, 1980.

(5) David Singer, *“Stable orbits and bifurcation of maps of the interval”* (SIAM vol. 35, No 2, September 1978)

(6) M. J. Feigenbaum, "Quantitative Universality for a Class of Non-Linear Transformations." *J. Stat. Phys.* **19**, 25-52, 1978.

(7) Andrei Vieru, ‘*Generalized iteration, Catastrophes, Generalized Sharkovsky’s Ordering’* arXiv:0801.3755 **[math.DS]**

**(**8) Andrei Vieru “Short note on additive sequences and on recursive processes” arxiv: 0806.3972